\documentclass{article}
\usepackage{lmodern}
\usepackage{graphicx} 
\usepackage[margin=1.5in]{geometry}
\usepackage{amsthm}
\usepackage{amssymb}
\usepackage{amsmath}
\usepackage{xcolor}
\usepackage{enumerate}
\usepackage{microtype}
\usepackage[numbers]{natbib}
\usepackage{hyperref}

\newtheorem{theorem}{Theorem}

\newtheorem{proposition}[theorem]{Proposition}
\newtheorem{corollary}[theorem]{Corollary}
\newtheorem{problem}{Problem}
\newtheorem{conjecture}{Conjecture}

\def\inst#1{$^{#1}$}

\title{A Survey on Ordered Ramsey Numbers\thanks{M. Balko was supported by grant no. 23-04949X of the Czech Science Foundation (GA\v{C}R) and by the Center for Foundations of Modern Computer Science (Charles Univ. project UNCE 24/SCI/008).
}}
\author{
Martin Balko\inst{1}
}

\begin{document}

\maketitle

\begin{center}
{\footnotesize
\inst{1} 
Department of Applied Mathematics, \\
Faculty of Mathematics and Physics, Charles University, Czech Republic \\
\texttt{balko@kam.mff.cuni.cz}
}
\end{center}

\begin{abstract}
The ordered Ramsey number of a graph $G^<$ with a linearly ordered vertex set is the smallest positive integer $N$ such that any
two-coloring of the edges of the ordered complete graph on $N$ vertices contains a monochromatic copy of $G^<$ in the given ordering.
The study of the quantitative behavior of ordered Ramsey numbers is a relatively new theme in Ramsey theory full of interesting and difficult problems.

In this survey paper, we summarize recent developments in the theory of ordered Ramsey numbers.
We point out connections to other areas of combinatorics and some well-known conjectures.
We also list several new and challenging open problems and highlight the often strikingly different behavior from the unordered case.

\end{abstract}

\section{Introduction}

In a broad sense, \emph{Ramsey theory} refers to any result whose underlying philosophy can be captured by the statement “Every structure of a given kind contains a large well-organized substructure”.
This part of discrete mathematics has developed enormously in the last few decades, emerging into a field with important statements in many areas, including combinatorics, geometry, logic, and number theory.

Our focus will be on graph theory where one of the cornerstones of Ramsey theory, from which Ramsey theory derives its name, is \emph{Ramsey's theorem}~\cite{ramsey29}.
This classical result states that for any two graphs $G$ and $H$, there is a positive integer $N$ such that every red-blue coloring of the complete graph $K_N$ on $N$ vertices contains a red copy of $G$ or a blue copy of $H$ as a subgraph of $K_N$.
The smallest such $N$ is called the \emph{Ramsey number} of $G$ and $H$ and is denoted by $R(G, H)$.
When $G$ and $H$ are isomorphic, we write $R(G)$ instead of $R(G, G)$.

Estimating the growth rate of Ramsey numbers is a notoriously difficult problem in general with a remarkable influence in combinatorics.
Despite several attempts over the last 70 years, see for example~\cite{conlon09,sah20,spencer75}, Ramsey numbers are not fully understood even for complete graphs.
For a long time, the best-known bounds on $R(K_n)$ essentially were 
\begin{equation}
\label{eq-ramseyBound}
2^{n/2} \leq R(K_n) \leq 2^{2n}.
\end{equation}
The lower bound is due to Erd\H{o}s~\cite{erdos47} while the upper bound was proved by Erd\H{o}s and Szekeres~\cite{erdosSzekeres35} in their seminal paper from 1935 where they independently rediscovered Ramsey's theorem.
Very recently, Campos, Griffiths, Morris,  and Sahasrabudhe~\cite{cgms23} made a breakthrough and improved the exponential upper bound to $R(K_n) \leq (4-\varepsilon)^n$ for some small constant $\varepsilon>0$.

For sparser graphs, Ramsey numbers grow at a significantly slower rate.
For example, a classical result by Chv\'{a}tal, R\"{o}dl, Szemer\'{e}di, and Trotter~\cite{crst83} says that Ramsey numbers of bounded-degree graphs are only linear in the number of vertices.
More precisely, for every positive integer $\Delta$, there is a positive integer $C=C(\Delta)$ such that every $n$-vertex graph $G$ with maximum degree $\Delta$ satisfies $R(G) \leq C \cdot n$.

In this survey, we focus on a more recent variant of Ramsey numbers defined for \emph{ordered graphs}, that is, graphs with linearly ordered vertex sets.
Studying these so-called \emph{ordered Ramsey numbers} became quite an active area in graph Ramsey theory~\cite{cfs15} which offers numerous applications in several parts of combinatorics such as discrete geometry and extremal theory of matrices.
We try to summarize the recent developments in the theory of ordered Ramsey numbers, pointing out connections to other areas of combinatorics and some well-known conjectures.
We also list several new and challenging open problems and highlight the often strikingly different behavior from the unordered case.

In Section~\ref{sec-orderedGraphs}, we survey the recent developments about ordered Ramsey numbers of ordered graphs.
We start by mentioning general bounds on ordered Ramsey numbers.
Then, we will focus on specific classes of ordered graphs, where we can prove much more precise estimates and, in some cases, even the exact formulas for the ordered Ramsey numbers.
We also discuss extensions of ordered Ramsey numbers for more than two colors.
In Section~\ref{sec-orderedHypergraphs}, we consider \emph{hypergraph ordered Ramsey numbers}.
Besides general bounds, we focus on so-called monotone paths, a case with close connections to the classical Erd\H{o}s--Szekeres theorem.
Finally, in Section~\ref{sec-edgeOrdered} we discuss \emph{edge-ordered Ramsey numbers} where we work with the ordering of the edges, not vertices.

We note that this paper is not intended to serve as an exhaustive survey of the subject.
For example, we treat problems of an asymptotic nature rather than being concerned with the computation of exact ordered Ramsey numbers.
We apologize in advance for any particularly glaring omissions.

For the sake of clarity of
presentation, we will maintain some conventions throughout the paper.
For a positive integer $n$, we use $[n]$ to denote the set $\{1,\dots,n\}$.
We omit floor and ceiling signs whenever they are not crucial and we use $\log$ to denote base $2$ logarithm.

\section{Graph Ordered Ramsey Numbers}
\label{sec-orderedGraphs}

An \emph{ordered graph} $G^<$ is a pair $(G,<)$ consisting of a graph $G$ and a linear ordering $<$ of the vertex set of $G$.
We call $G^<$ an \emph{ordering} of $G$.
Two ordered graphs $G^<$ and $H^<$ are \emph{isomorphic} if their underlying graphs $G$ and $H$ are isomorphic via a one-to-one correspondence $V(G) \to V(H)$ that preserves the vertex orderings.
Note that, for every positive integer $n$, there is only one ordered complete graph on $n$ vertices up to isomorphism.
We use $K^<_n$ to denote such an ordered complete graph.

An ordered graph $G^<$ is an \emph{ordered subgraph} of an ordered graph $H^<$
if $G$ is a subgraph of $H$ and the vertex ordering of $G$ is a suborder of the vertex ordering of $H$. 
If an ordered graph $H^<$ contains an ordered subgraph isomorphic to an ordered graph $G^<$, then we say that $H^<$ contains a \emph{copy} of $G^<$.

We can now define an analog of Ramsey numbers for ordered graphs.
The \emph{ordered Ramsey number} $R_<(G^<, H^<)$ of two ordered graphs $G^<$ and $H^<$ is the minimum positive integer $N$ such that every red-blue coloring of the edges of~$K^<_N$ contains a red copy of $G^<$ or a blue copy of~$H^<$.
If $G^<$ and $H^<$ are isomorphic, then we simply write $R_<(G^<)$ instead of $R_<(G^<,H^<)$ and we call it the \emph{diagonal} case.
We refer to $R_<(G^<,H^<)$ with non-isomorphic $G^<$ and $H^<$ as the \emph{off-diagonal} case.

Note that we have $R(G)\leq R_<(G^<)$ for every graph $G$ and each ordering~$G^<$ of $G$.
Moreover, $R(K_n)=R_<(K^<_n)$ for every positive integer $n$.
Since every ordered graph $G_<$ on $n$ vertices is an ordered subgraph of~$K^<_n$, we thus obtain the bounds 
\[R(G) \leq R_<(G^<)\leq R(K_n).\]
It follows that the number $R_<(G^<)$ is finite for every ordered graph $G^<$.
Similarly, $R_<(G^<, H^<)$ is also finite for any ordered graphs $G^<$ and $H^<$ and thus ordered Ramsey numbers are well-defined.

\begin{figure}[htb]
\centering
\includegraphics{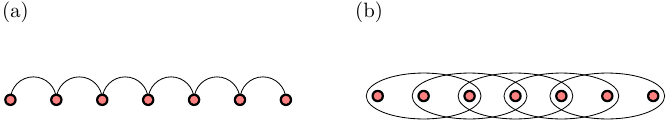} 
\caption{(a) The monotone path $MP^<_n$ for $n=7$. (b) The monotone path $MP^{<(3)}_n$ for $n=7$.
In all figures of ordered graphs, the vertices are ordered from left to right.}
\label{fig-monPath}
\end{figure}

From some point of view, the study of ordered Ramsey numbers is as old as Ramsey theory itself.
For example, one of the oldest Ramsey-type results, the famous \emph{Erd\H{o}s--Szekeres lemma} on monotone subsequences proved by Erd\H{o}s and Szekeres~\cite{erdosSzekeres35} in 1935, can be derived from estimates on ordered Ramsey numbers. 
This classical result states that for every positive integer $n$, every sequence of at least $(n-1)^2+1$ distinct real numbers contains an increasing or a decreasing subsequence of length $n$.
Moreover, this bound is sharp.
To see that it follows from results in ordered Ramsey theory, consider the following ordering $MP^<_n$, called the \emph{monotone path}, of the path $P_n$ on $n$ vertices.
If $v_1 < \dots < v_n$ are the vertices of~$MP^<_n$, then the edges of $MP^<_n$ are the pairs $\{v_i,v_{i+1}\}$ for every $i=1,\dots,n-1$; see part~(a) of Figure~\ref{fig-monPath}.
Given a sequence $S=(s_1,\dots,s_N)$ of distinct real numbers, we construct an ordered graph $K^<_n$ with vertex set $S$ and the ordering of the vertices given by their positions in~$S$.
Then, we color an edge $\{s_i,s_j\}$ with $i<j$ red if $s_i < s_j$ and blue otherwise.
Afterward, red monotone paths on $n$ vertices correspond to increasing subsequences of $S$ of length $n$ and blue monotone paths on $n$ vertices to decreasing subsequences of $S$ of length $n$.
The Erd\H{o}s--Szekeres lemma now follows from the fact 
\begin{equation}
\label{eq-ordRam-monotonePaths}
R_<(MP^<_n) = (n-1)^2+1
\end{equation}
proved, for example, by Choudum and Ponnusamy~\cite{choudum02} or by Milans, Stolee, and West~\cite{milans12}.

Similarly, the Erd\H{o}s--Szekeres theorem on convex polygons can be derived from estimates on ordered Ramsey numbers of 3-uniform monotone hyperpaths; see Subsection~\ref{subsec-monotonePaths}.
The ordered Ramsey numbers are also closely connected to the extremal theory of $\{0,1\}$-matrices~\cite{bckk13,jjmm24} and variants of Ramsey numbers from discrete geometry, for example, to so-called \emph{geometric Ramsey numbers}~\cite{kar97,kar98}.
Given these results and connections, there is a strong motivation to study ordered Ramsey numbers and their variants.

While there had been a lot of results about ordered Ramsey numbers of monotone hyperpaths~\cite{choudum02,elias11,fox12,milans12,moshkovitz12}, there was surprisingly little work on ordered Ramsey numbers of more general ordered graphs and hypergraphs. 
The first systematic study of ordered Ramsey numbers was conducted by Balko, Cibulka, Kr\'{a}l, and Kyn\v{c}l~\cite{bckk13} and independently by Conlon, Fox, Lee, and Sudakov~\cite{clfs17}.
Since then there has been much progress in understanding the ordered Ramsey numbers, which we describe in this section.

\subsection{Bounded degrees and interval chromatic number}
\label{subsec-boundedDegInt}

Since $R_<(G^<) \leq R(K_n)$ for every ordered graph $G^<$ on $n$ vertices, we see from~\eqref{eq-ramseyBound} that $R_<(G^<)$ grows at most exponentially in $n$.
For dense ordered graphs, a standard probabilistic argument shows that this is asymptotically tight.
However, the ordered Ramsey numbers may differ substantially from the usual Ramsey numbers for sparser graphs.
This result was proved independently by Balko, Cibulka, Kr\'{a}l, and Kyn\v{c}l~\cite{bckk13} and by Conlon, Fox, Lee, and Sudakov~\cite{clfs17} who showed that there are ordered matchings with superpolynomial ordered Ramsey numbers.
A \emph{matching} is a graph consisting of pairwise disjoint edges.

\begin{theorem}[\cite{bckk13,clfs17}]
\label{thm-oderRam-match}
There is a constant $C>0$ such that for every even $n \geq 2$, there is an ordered matching $M^<_n$ on $n$ vertices satisfying
\[R_<(M^<_n) \geq n^{C\log{n}/\log{\log{n}}}.\]
\end{theorem}

Note that this result is in sharp contrast with the classical result by Chv\'{a}tal, R\"{o}dl, Szemer\'{e}di, and Trotter~\cite{crst83}, according to which Ramsey numbers of bounded-degree graphs are only linear.
For ordered matching, this is a rather typical behavior, as Conlon, Fox, Lee, and Sudakov~\cite{clfs17} proved that this superpolynomial lower bound holds for almost all ordered matchings.
Very recently, Fox and He~\cite{foxHe23} applied Theorem~\ref{thm-oderRam-match} to obtain superpolynomial lower bounds on Ramsey numbers of sparse digraphs.

Given the lower bound from Theorem~\ref{thm-oderRam-match}, it is clear that to obtain polynomial upper bounds on ordered Ramsey numbers, bounding the maximum degree is not enough and we need to use some other graph parameter.
Balko, Cibulka, Kr\'{a}l, and Kyn\v{c}l~\cite{bckk13} and independently Conlon, Fox, Lee, and Sudakov~\cite{clfs17} showed that if we additionally bound the following analog of the chromatic number, then we obtain polynomial upper bounds on ordered Ramsey numbers.
We will state the bound obtained in~\cite{clfs17}, as it is stronger than the one in~\cite{bckk13}, but first, we need to introduce some definitions.

A subset $I$ of vertices of an ordered graph $G^<$ is an \emph{interval} in $G^<$ if for every pair $u,v$ of vertices of $I$ with $u < v$, every vertex $w$ of $G^<$ satisfying $u < w < v$ is contained in~$I$.
The \emph{interval chromatic number} of $G^<$ is the minimum number of intervals the vertex set of $G^<$ can be partitioned into so that there is no edge between vertices of the same interval.
We note that there is a variant of the \emph{Erd\H{o}s--Stone--Simonovits theorem} for ordered graphs proved by Pach and Tardos~\cite{pachTardos06}, which is expressed in terms of the interval chromatic number.
For a positive integer $d$, a graph $G$ is \emph{$d$-degenerate} if every subgraph of $G$ contains a vertex of degree at most $d$.
For positive integers $t,n_1,\dots,n_t$, let $K^<_{n_1,\dots,n_t}$ be an ordering of the complete $t$-partite graph $K_{n_1,\dots,n_t}$ in which the vertices of the color class of size $n_i$ form the $i$th interval.
If $n_1=\dots=n_t=n$, we simply write $K^<_t(n)$ instead of $K^<_{n_1,\dots,n_t}$.

\begin{theorem}[\cite{clfs17}]
\label{thm-oderRam-degreeIntNum2}
Let $G^<$ be an ordered $d$-degenerate graph on $n$ vertices with maximum degree $\Delta$.
For positive integers $n'$ and $\chi$, let $s=\lceil\log{\chi} \rceil$ and $D = 8\chi^2n'$.
Then,
\[R_<(G^<,K^<_{\chi}(n')) \leq 2^{s^2d+s}\Delta^sn^sD^{ds+1}.\]
\end{theorem}

In particular, if $G^<$ is an ordered $d$-degenerate graph with $n$ vertices and with interval chromatic number $\chi$, then Theorem~\ref{thm-oderRam-degreeIntNum2} implies
\begin{equation}
\label{eq-ordRam-degInt}
R_<(G^<) \leq n^{32d\log{\chi}}.
\end{equation}
Note that for fixed $d$ this bound almost matches the lower bound from Theorem~\ref{thm-oderRam-match}.
Although the bounds are quite close, there is still a small gap in the exponent, even for $d=1$.
Conlon, Fox, Lee, and Sudakov~\cite{clfs17} thus posed the following problem.

\begin{problem}[\cite{clfs17}]
Close the gap between the upper and lower bounds for ordered Ramsey numbers of
matchings.
\end{problem}

It follows from~\eqref{eq-ordRam-degInt}  that bounded-degree ordered graphs with bounded interval chromatic numbers have polynomial ordered Ramsey numbers.
However, we have no non-trivial lower bounds for this case and thus Balko, Cibulka, Kr\'{a}l, and Kyn\v{c}l~\cite{bckk13} stated the following problem.

\begin{problem}[\cite{bckk13}]
Is there a constant $c>0$ such that for every fixed $\Delta$ there is a sequence $\{G^<_n\}_{n\in\mathbb{N}}$ of ordered $\Delta$-regular graphs $G^<_n$ with $n$ vertices and interval chromatic number $2$ such that $R_<(G^<_n) \ge n^{c\Delta}$?
\end{problem}

Using a result of Erd\H{o}s and Szemer\'{e}di~\cite{erdosSzemeredi72}, Conlon, Fox, Lee, and Sudakov~\cite{clfs17} derived the following result from the proof of Theorem~\ref{thm-oderRam-degreeIntNum2}, which shows that the ordered Ramsey numbers behave more like the usual Ramsey numbers for dense ordered graphs.

\begin{theorem}[\cite{clfs17}]
\label{thm-oderRam-denser}
There is a constant $C>0$ such that every ordered $d$-degenerate graph $G^<$ on $n$ vertices satisfies
\[R_<(G^<) \leq 2^{Cd\log^2{(2n/d)}}.\]
\end{theorem}

This result is close to sharp for very small $d$ by Theorem~\ref{thm-oderRam-match} and for very large $d$ by~\eqref{eq-ramseyBound}.

A simple argument shows that if an $n$-vertex ordered matching $M^<$ has the interval chromatic number 2, then the bound $R_<(M^<) \leq n^{\lceil \log{n}\rceil}$ proved by Conlon, Fox, Lee, and Sudakov~\cite{clfs17} for any $n$-vertex ordered matching $M^<$ can be significantly improved to $O(n^2)$~\cite{clfs17}.
This is quite close to the truth, as Conlon, Fox, Lee, and Sudakov~\cite{clfs17} constructed an ordered matching $M^<$ on $n$ vertices with interval chromatic number 2 such that $R_<(M^<) \geq \frac{cn^2}{\log^2{(n)}\log{\log{n}}}$ for some constant $c>0$.
This was improved by Balko, Jel\'{i}nek, and Valtr~\cite{bjv16}, who proved a stronger bound, which additionally holds for almost all ordered matchings with interval chromatic number 2.
To state this result, we need to introduce some definitions first, which will also be useful later.

For a positive integer $n$, the \emph{random $n$-permutation} is a permutation of the set $[n]$ chosen independently uniformly at random from the set of all $n!$ permutations of the set $[n]$.
For a positive integer $n$ and the random $n$-permutation $\pi$, the \emph{random ordered matching $M^<(\pi)$} on $n$ vertices is the ordered matching with the vertex set $[2n]$ and with edges $\{i,n+\pi(i)\}$ for every $i \in [n]$.
Note that the interval chromatic number of every random ordered matching is~2.
The random ordered matching satisfies an event $A$ \emph{asymptotically almost surely} if the probability that $A$ holds tends to 1 as $n$ goes to infinity.

\begin{theorem}[\cite{bjv16}]
\label{thm-oderRam-bipMatching}
There is a constant $C>0$ such  that the random ordered matching $M^<(\pi)$ on $n$ vertices asymptotically almost surely satisfies 
\[R_<(M^<(\pi)) \geq C \cdot \left(\frac{n}{\log{n}}\right)^2.\]
\end{theorem}

Similarly, as in the case of general ordered matchings, there is a gap between the lower and upper bound.
Thus, Colon, Fox, Lee, and Sudakov~\cite{clfs17} posed the following problem for ordered matchings with interval chromatic number 2.

\begin{problem}[\cite{clfs17}]
Close the gap between the upper and lower bounds for ordered Ramsey numbers of ordered
matchings with interval chromatic number 2.
\end{problem}

Geneson, Holmes, Liu, Neeidinger, Pehova, and Wass~\cite{ghnpw19} posed a similar problem for ordered paths with interval chromatic number 2.
Here, it follows from Theorem~\ref{thm-oderRam-degreeIntNum2} that the ordered Ramsey number $R_<(P^<)$ of any ordered path $P^<$ on $n$ vertices with interval chromatic number 2 is at most $O(n^3)$ while the best known lower bound $\Omega((n/\log{n})^2)$ follows from Theorem~\ref{thm-oderRam-bipMatching}.

\begin{problem}[\cite{ghnpw19}]
Is it true that $R_<(P^<) \leq O(n^2)$ for every ordering $P^<$ of the path on $n$ vertices with interval chromatic number 2?
\end{problem}

\subsection{Bounded bandwidth}

Besides the interval chromatic number, another natural parameter for ordered graphs is their bandwidth.
For an ordered graph $G^<$, the \emph{bandwidth} of $G^<$ is the length of the longest edge in $G^<$.
That is, the maximum of $|i-j|$ taken over all edges $\{u,v\}$ of $G^<$, where $i$ is the position of $u$ and $j$ is the position of~$v$ in~$<$.
We call the number $|i-j|$ the \emph{length} of the edge $\{u,v\}$.

Conlon, Fox, Lee, and Sudakov~\cite{clfs17} proved that, for every positive integer $k$, every ordered matching $M^<$ on $n$ vertices with bandwidth at most $k$ satisfies $R_<(M^<) \leq n^{\lceil \log{k} \rceil+2}$.
They also asked whether this result can be extended by proving a polynomial upper bound on ordered Ramsey numbers of all ordered graphs with bounded bandwidth.
This problem was solved by Balko, Cibulka, Kr\'{a}l, and Kyn\v{c}l~\cite{bckk13} who proved that, for every positive  integer $k$,  there is a constant $C=C(k)$ such  that  every $n$-vertex ordered graph $G^<$ with bandwidth $k$ satisfies
$R_<(G^<) \leq C \cdot n^{128k}$.

Observe that every $n$-vertex ordered graph $G^<$ with bandwidth at most $k$ is an ordered subgraph of the $n$-vertex ordered graph $P^<_{k,n}$ that contains all edges of length at most $k$.
In particular, $R_<(G^<) \leq R_<(P^<_{k,n})$.
Note that $P^<_{1,n} = MP^<_n$ and thus $R_<(P^<_{1,n}) = (n-1)^2+1$ by~\eqref{eq-ordRam-monotonePaths}.
Mubayi~\cite{mubayi17} improved the upper bound on $R_<(P^<_{k,n})$ in the case $k=2$ by showing $R_<(P^<_{2,n}) \leq O(n^{19.487})$ for every $n \geq 2$.
He also used this result to determine the correct tower growth rate of the $k$-uniform hypergraph Ramsey number of a $(k+1)$-clique versus a monotone $k$-uniform path; see Section~\ref{sec-orderedHypergraphs} for the definitions.
After that, Gishboliner, Jin, and Sudakov~\cite{gishJinSud24} proved that $R_<(P^<_{k,n}) \leq c n^{4k-2}$ for some constant $c=c(k)>0$ and asked whether in fact $R_<(P^<_{k,n}) \leq c' n^C$ for a constant $c'=c'(k)>0$ and an absolute constant $C>0$.
Very recently, Gir\~{a}o, Janzer, and Janzer~\cite{gjj24} gave a positive answer to this problem and derived the currently strongest bound.

\begin{theorem}[\cite{gjj24}]
For every $k \in \mathbb{N}$ and $\varepsilon >0$, there is $C=C(k,\varepsilon)>0$ such that for every positive integer $n$ we have
\[R_<(P^<_{k,n}) \leq C \cdot n^{4+\varepsilon}.\]
\end{theorem}

Motivated by Mubayi's work~\cite{mubayi17} on hypergraph Ramsey numbers, Mubayi and Suk~\cite{mubSuk24} considered the off-diagonal setting and proved that $R_<(P^<_{k,n},K^<_n) \leq (2n)^{k(k+1)\log{n}}$ for all positive integers $k$ and $n$.
 Gishboliner, Jin and Sudakov~\cite{gishJinSud24} proved that $R_<(P^<_{k,n},K^<_n) \leq c n^{k(2k-1)}$ for some constant $c=c(t)>0$ and conjectured that the exponent can
even be chosen to be linear in $k$.
This conjecture was verified by Gir\~{a}o, Janzer, and Janzer~\cite{gjj24} who proved  that there exists $C > 0$ such that, for all positive integers $n > k$,
\[R_<(P^<_{k,n},K^<_n) \leq n^{Ck}.\]
In fact, Gir\~{a}o, Janzer, and Janzer~\cite{gjj24} derived an essentially tight bound
\[R_<(P^<_{k,s},K^<_n) \leq R(K_s,K_k)^C \cdot n\]
which is valid
for all $s$, $t$ and $n$, improving earlier estimate by Mubayi and Suk~\cite{mubSuk24}. 

Despite this remarkable progress, there are still some open problems left.
For example, there is still a gap between the lower and upper bounds on $R_<(P^<_{k,n})$.
While the best upper bound is of order $n^{4+o(1)}$, we only have a trivial quadratic lower bound.
Gir\~{a}o, Janzer, and Janzer~\cite{gjj24} believed that the truth is closer to the lower bound and posed the following problem.

\begin{problem}[\cite{gjj24}]
Is $R_<(P^<_{k,n})= n^{2+o(1)}$ as $n \to \infty$? 
\end{problem}

\subsection{Off-diagonal Ordered Ramsey Numbers}

For the off-diagonal ordered Ramsey numbers, Conlon, Fox, Lee, and Sudakov~\cite{clfs17} investigated the ordered Ramsey numbers $R_<(M^<, K^<_3)$, where $M^<$ is an ordered matching on $n$ vertices.
It follows from the well-known bound $R(K_n,K_3) \leq O(n^2/\log{n})$ by Ajtai, Koml\'{o}s, and Szemer\'{e}di~\cite{ajKomSze80} that \[R_<(M^<,K^<_3) \leq O\left(\frac{n^2}{\log{n}}\right).\]
For $R(K_n,K_3)$, this bound is tight as shown by Kim~\cite{kim95}, but Conlon, Fox, Lee, and Sudakov~\cite{clfs17} expect that
this upper bound is far from optimal for $R_<(M^<,K^<_3)$.
They constructed an ordered matching $M^<$ on $n$ vertices satisfying \[R_<(M^<,K^<_3) \geq O\left(\left(\frac{n}{\log{n}}\right)^{4/3}\right)\] 
and posed the following problem.

\begin{problem}[\cite{clfs17}]
\label{prob-oderRam-matchingK3}
Does there exist an $\varepsilon > 0$ such that every ordered matching $M^<$ on $n$ vertices satisfies 
$R_<(M^<,K^<_3) \leq O(n^{2-\varepsilon})$?
\end{problem}

In case there is an affirmative answer to this problem, one has to employ the sparsity of ordered matching to distinguish it from $K_n$.
However, ordered matchings on $n$ vertices typically contain an edge between any two intervals of length $3\sqrt{n}\log{n}$~\cite{clfs17} and thus can be used to ``simulate'' complete graphs.
This makes Problem~\ref{prob-oderRam-matchingK3} particularly interesting.
However, it is still open and seems to be difficult.

Nevertheless, there has been some partial progress.
Rohatgi~\cite{rohatgi19} proved that asymptotically almost surely the random ordered matching $M^<(\pi)$ on $n$ vertices with interval chromatic number~2 satisfies $R_<(M^<(\pi),K^<_3) \leq O(n^{24/13})$.
He also showed that if $M^<$ is a non-crossing ordered matching, then, for every $\varepsilon > 0$, we have $R_<(M^<,K^<_3) \leq cn^{1+\varepsilon}$ for some constant $c=c(\varepsilon)> 0$.
Here, an ordered graph $G^<$ is \emph{non-crossing} if it does not contain two edges $\{u,v\}$ and $\{x,y\}$ with $u < x < v < y$.
Balko and Poljak~\cite{balPol24} later showed that $R_<(M^<(\pi),K^<_3) \geq \Omega((n/\log n)^{5/4})$ and $R_<(M^<(\pi),K^<_3) \leq O(n^{7/4})$ asymptotically almost surely, improving the result by Rohatgi~\cite{rohatgi19}.

Since Problem~\ref{prob-oderRam-matchingK3} seems difficult, some weaker variants of this problem were introduced in the literature.
For example, Rohatgi~\cite{rohatgi19} posed the following interesting conjecture.

\begin{conjecture}[\cite{rohatgi19}]
For a positive integer $\chi$, there is a constant $\varepsilon(\chi)>0$ such that 
\[R_<(M^<,K^<_3) \leq O(n^{2-\varepsilon(\chi)})\]
for almost every ordered matching $M^<$ on $n$ vertices with interval chromatic number~$\chi$.
\end{conjecture}

Even in the case of $M^<(\pi)$ the bounds are roughly by the factor $\sqrt{n}$ apart, so Balko and Poljak~\cite{balPol24} asked the following question.

\begin{problem}[\cite{balPol24}]
\label{prob-random}
What is the growth rate of $R_<(M^<(\pi),K^<_3)$ for the random ordered matchings $M^<(\pi)$ on $n$ vertices with interval chromatic number 2?
\end{problem}

\subsection{Minimum Ordered Ramsey Numbers}

Conlon, Fox, Lee, and Sudakov~\cite{clfs17} characterized graphs~$G$ for which the ordered Ramsey number $R_<(G^<)$ is linear in the number of vertices of $G$ for every ordering $G^<$ of $G$.
These are precisely graphs~$G$ whose edges can be covered by a constant number of vertices.
A similar problem is to determine graphs that admit an ordering for which the corresponding ordered Ramsey number is linear.
This motivated the following definition.

For an unordered graph $G$, the \emph{minimum ordered Ramsey number of $G$} is defined as
\[\min R_<(G)=\min\{R_<(G^<): G^<\ \text{is an ordering of}\ G\}.\]

Since Ramsey numbers of bounded-degree graphs are linear in the number of vertices, it is natural to ask whether the minimum ordered Ramsey numbers of bounded-degree graphs are always at most linear.
Conlon, Fox, Lee, and Sudakov~\cite{clfs17} considered this unlikely and asked whether random 3-regular graphs have superlinear ordered Ramsey numbers for all orderings.

Balko, Jel\'{i}nek, and Valtr~\cite{bjv16} gave an affirmative answer to this problem. 
For a positive integer $d$, we let $G(d,n)$ denote the random $d$-regular graph on $n$ vertices drawn uniformly and independently from the set of all $d$-regular graphs on the vertex set $[n]$.

\begin{theorem}[\cite{bjv16}]
\label{thm-ordRam-lowerBoundAllOrderings}
For every fixed integer $d \geq 3$, asymptotically almost surely
\[\min R_<(G(d,n)) \geq \frac{n^{3/2-1/d}}{4\log{n}\log{\log{n}}}.\]
\end{theorem}

Theorem~\ref{thm-ordRam-lowerBoundAllOrderings} shows that random $d$-regular graphs have superlinear minimum ordered numbers for any fixed integer number $d \geq 3$.
In particular, almost every 3-regular graph $G$ on $n$ vertices satisfies 
\[\min R_<(G) \geq n^{7/6}/(4\log{n}\log{\log{n}}).\]

It can be easily seen that the minimum ordered Ramsey numbers of ordered matchings are linear.
Balko, Jel\'{i}nek, and Valtr~\cite{bjv16} showed that this is the case for 2-regular graphs as well by proving that there is a constant $C$ such that for every graph~$G$ on $n$ vertices with maximum degree 2, we have $\min R_<(G) \leq Cn$.
In fact, Balko, Jel\'{i}nek, and Valtr~\cite{bjv16} proved a stronger Tur\'{a}n-type statement for bipartite graphs of maximum degree 2.

For the upper bounds in the case of a larger maximum degree, a simple corollary of Theorem~\ref{thm-oderRam-degreeIntNum2} states that every graph $G$ on $n$ vertices with constant maximum degree $\Delta$ admits an ordering $G^<$ with $R_<(G^<)$ polynomial in~$n$.
More precisely, every graph $G$ with $n$ vertices and with maximum degree $\Delta$ satisfies \begin{equation}
\label{eq-ordRam-minRamseyUpper}
\min R_<(G) \leq O(n^{(\Delta+1)\lceil\log(\Delta+1)\rceil+1}).
\end{equation}

Thus, the currently best-known bounds on minimum ordered Ramsey numbers of 3-regular graphs are of order $\Omega(n^{7/6}/(\log{n}\log{\log{n}}))$ and $O(n^{4\lceil\log(4)\rceil+1}) = O(n^9)$ by Theorem~\ref{thm-ordRam-lowerBoundAllOrderings} and by~\eqref{eq-ordRam-minRamseyUpper}, respectively.
Since the gap is rather large, Balko, Jel\'{i}nek, and Valtr~\cite{bjv16} posed the following problem.

\begin{problem}[\cite{bjv16}]
Close the gap between the upper and lower bounds for minimum ordered Ramsey numbers of 3-regular graphs.
\end{problem}

\subsection{Specific Classes of Ordered Graphs}

Obtaining precise formulas for ordered Ramsey numbers is difficult, essentially the only non-trivial class of graphs with fully characterized ordered Ramsey numbers are stars~\cite{bckk13,choudum02}.
However, there are some specific and natural orderings for which it is possible to derive such formulas or at least very close bounds that often lead to interesting applications.
In this section, we focus on three classes of ordered graphs: nested matchings, alternating paths, and monotone cycles.

\subsubsection{Nested Matchings}

\begin{figure}[htb]
\centering
\includegraphics{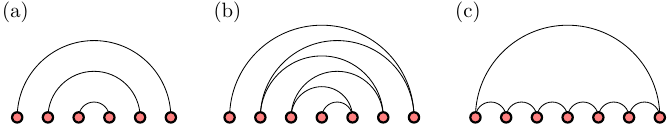} 
\caption{(a) The nested matching $NM_k$ for $k=3$. (b) The alternating path $AP^<_n$ for $n=7$. (c) The monotone cycle $MC^<_n$ for $n=7$.}
\label{fig-specificClasses}
\end{figure}

A basic building block in the proof of the result by Rohatgi~\cite{rohatgi19} about the ordered Ramsey numbers $R_<(M^<, K^<_3)$ for a non-crossing ordered matching $M^<$ is based on so-called nested matchings.
A \emph{nested ordered matching} $NM^<_k$ on vertices $v_1 < \dots < v_{2k}$ has edges $\{v_i,v_{n-i+1}\}$ for every $i=1,\dots,k$; see part~(a) of Figure~\ref{fig-specificClasses}.

Rohatgi~\cite{rohatgi19} proved 
$4k-1 \leq R_<(NM^<_k,K^<_3) \leq 6k$
for every positive integer~$k$.
He believed that the lower bound was tight and conjectured that $R_<(NM^<_k, K^<_3) = 4k-1$.
This conjecture is true for $k \leq 3$, but Balko and Poljak~\cite{balkoPoljak21} disproved it for any $k \geq 4$ by showing the following bounds.

\begin{theorem}[\cite{balkoPoljak21}]
\label{thm-oderRam-nestedMatchings}
For every positive integer $k$, we have 
\[R_<(NM^<_k,K^<_3) \leq (3+\sqrt{5})k < 5.3k.\]
If $k \geq 6$, we have
\[R_<(NM^<_k,K^<_3) \geq 4k+1.\] 
\end{theorem}

Using the lower bounds from Theorem~\ref{thm-oderRam-nestedMatchings}, Balko and Poljak improved the best-known bounds on the maximum chromatic number of so-called \emph{$k$-queue graphs}, which addresses a problem posed by Dujmovi\`{c} and Wood~\cite{dujWoo10}.
Despite this progress, no precise formula for $R_<(NM_k, K^<_3)$ is known.

\subsubsection{Alternating paths}

The \emph{alternating path} $AP^<_n$ is the ordering of $P_n$ on vertices $v_1 < \dots < v_n$ with edges $\{v_i,v_j\}$ for $i+j\in \{n+1,n+2\}$; see part~(b) of Figure~\ref{fig-specificClasses}.
Note that the alternating path $AP^<_n$ is an ordered subgraph of $K^<_{\lceil n/2 \rceil,\lfloor n/2 \rfloor}$, and so it has interval chromatic number $2$. 

Using a result from the extremal theory of $\{0,1\}$-matrices, Balko, Cibulka, Kr\'{a}l, and Kyn\v{c}l~\cite{bckk13} proved the following linear bounds on~$R_<(AP^<_n)$.

\begin{proposition}[\cite{bckk13}]
\label{prop-ordRam-alt_paths}
For every integer $n>2$, we have 
\[5\lfloor n/2 \rfloor -4 \le R_<(AP^<_n) \le 2n-3+\sqrt{2n^2-8n+11}.\]
\end{proposition}

The proof of the lower bound was later extended by Neidinger and West~\cite{neiWes19} to obtain various linear lower bounds for various classes of ordered graphs with interval chromatic number 2.

The precise multiplicative factor in $R_<(P^<_n)$ is unknown.
Computer experiments by Balko, Cibulka, Kr\'{a}l, and Kyn\v{c}l~\cite{bckk13} indicate that ordered Ramsey numbers of alternating paths might be minimal among all path orderings, making it a particularly interesting class of ordered paths.

\begin{problem}[\cite{bckk13}]
\label{prob-ordRam-altMin}
For some positive integer $n$, is there an ordering $P^<_n$ of the path $P_n$ on $n$ vertices such that $R_<(P^<_n) < R_<(AP^<_n)$?
\end{problem} 

From the cases with $n \leq 8$, it seemed that $R_<(P^<_n)$ could be equal to $\lfloor (n-2)\frac{1+\sqrt{5}}{2}\rfloor+n$, but this was later disproved by Balko and Poljak~\cite{balPol24}.
We also note that the alternating path $AP^<_n$ is contained in every ordered graph on $N\geq n/\varepsilon$ vertices with at least $\varepsilon N^2$ edges~\cite{bjv16}.

\subsubsection{Monotone Cycles}

We do not know any precise formulas for ordered Ramsey numbers of ordered cycles except one particularly natural ordering of the cycle $C_n$, called the \emph{monotone cycle} $MC^<_n$, which is obtained from the monotone path $MP^<_n$ by adding an edge between the first and the last vertex in the vertex ordering; see part~(c) of Figure~\ref{fig-specificClasses}.

The precise formula for ordered cycles covers even the off-diagonal case and is given by the following result by Balko, Cibulka, Kr\'{a}l, and Kyn\v{c}l~\cite{bckk13}.

\begin{theorem}[\cite{bckk13}]
\label{thm-ordRam-monCycle}
For all integers $r \ge 2$ and $s \ge 2$, we have \[R_<(MC^<_r,MC^<_s)=2rs-3r-3s+6.\]
\end{theorem}

Theorem~\ref{thm-ordRam-monCycle} has an interesting application in discrete geometry.
The \emph{geometric Ramsey numbers}~\cite{cibGao13,kar97, kar98} are natural analogs of ordered Ramsey numbers.
For a finite set of points $P \subset \mathbb{R}^2$ in general position, let $K_P$ be the \emph{complete geometric graph on $P$}, which is a complete graph drawn in the plane so that its vertices are represented by the points in $P$ and the edges are drawn as straight-line segments between the pairs of points in $P$.
The graph $K_P$ is \emph{convex} if $P$ is in convex position.
The \emph{geometric Ramsey number} of a graph~$G$, denoted by $Rg(G)$, is the smallest $N$ such that every complete geometric graph $K_P$ on $N$ vertices with edges colored by two colors contains a non-crossing monochromatic drawing of~$G$.
If we consider only convex complete geometric graphs $K_P$ in the definition, then we get so-called \emph{convex geometric Ramsey number} $Rc(G)$.
Note that these numbers are finite only if $G$ is outerplanar and that $Rc(G) \le Rg(G)$ for every outerplanar graph $G$.

Balko, Cibulka, Kr\'{a}l, and Kyn\v{c}l~\cite{bckk13} observed that the geometric and convex geometric Ramsey numbers of cycles are equal to the ordered Ramsey numbers of monotone cycles, which together with Theorem~\ref{thm-ordRam-monCycle} gives the following formula for $Rc(C_n)$ and $Rg(C_n)$.

\begin{corollary}[\cite{bckk13}]
\label{cor-ordRam-monCycl}
For every integer $n \ge 3$, we have $Rc(C_n)=Rg(C_n)=2n^2 - 6n + 6$.
\end{corollary}

So in this case, vertex orderings are enough to determine these numbers, no deeper geometric structure is needed, which sometimes seems to be the case; see Subsection~\ref{subsec-monotonePaths} for other examples.

\subsection{Multicolor Ordered Ramsey Numbers}

The definition of Ramsey numbers extends to a multicolor case.
For any integer $q \geq 2$ and a graph $G$, we can define the \emph{$q$-color Ramsey number $R(G;q)$} as the smallest positive integer $N$ such that every $q$-coloring of $K_N$  contains a monochromatic copy of $G$ as a subgraph of $K_N$.
A straightforward modification of the Erdős and Szekeres neighborhood-chasing argument~\cite{erdosSzekeres35} yields
$R(K_n;q) \leq q^{qn}$.
On the other hand, a
product argument by Lefmann~\cite{lefmann1987lower} shows that $R(K_n;q) \geq 2^{qn/4}$.
Although there has recently been remarkable progress in obtaining better lower and upper bounds on $R(K_n;q)$~\cite{bbcghmst24,conlon2021multicolorlower,sawin2022lower,wigderson2021lower}, the $q$-color Ramsey numbers are still far from being completely understood.

The multicolor setting can be introduced for ordered Ramsey numbers as well.
For any integer $q \geq 2$ and ordered graphs $G^<_1,\dots,G^<_q$, the \emph{$q$-color ordered Ramsey number $R_<(G^<_1,\dots,G^<_q;q)$} is the smallest positive integer $N$ such that for every $q$-coloring of $K^<_N$ by colors from $[q]$ there is $i \in [q]$ and a copy of $G_i$ in color $i$ forming a ordered subgraph of $K^<_N$.
In the case $G^< = G^<_1 = \cdots = G^<_q$, we simply write $R_<(G^<;q)$ instead of $R_<(G^<_1,\dots,G^<_q;q)$.

Similarly as before, it can be shown that the numbers $R_<(G^<_1,\dots, G^<_q;q)$ are always bounded and that $R_<(G^<;q) \leq R(K_n;q)$ for any $n$-vertex ordered graph~$G^<$.
However, we only have a poor understanding of ordered Ramsey numbers for more than two
colors, even for ordered matchings.
Conlon, Fox, Lee, and Sudakov~\cite{clfs17} proved the upper bound $R_<(M^<; q) \leq  n^{(2 \log n)^{(q-1)}}$ for any ordered matching $M^<$ on $n$ vertices.
They believe that a much stronger upper bound should hold and pose the following problem.

\begin{problem}[\cite{clfs17}] For any integer $q \geq 3$, show that there exists a constant $c_q$ such that $R_<(M^<; q) \leq n^{c_q \log n}$ for any ordered matching $M^<$ on $n$ vertices.
\end{problem}

For ordered matchings on $n$ vertices with interval chromatic number two and $q \geq 2$ colors, Balko and Grinerov\'{a}~\cite{balGrin24} showed that $R_<(M^<; q) = n^{\Theta (q)}$; see also~\cite{grin24}. 
It would be interesting to determine the leading constant in the exponent. 

\begin{problem}[\cite{balGrin24,grin24}] For any integer $q \geq 3$, determine the constant $c$ such that $R_<(M^<; q) = n^{c q + o(1)}$ for any ordered matching $M^<$ on $n$ vertices with interval chromatic number two.
\end{problem}

In the off-diagonal setting, Balko and Grinerov\'{a}~\cite{balGrin24,grin24} also introduced a multicolor variant of Problem~\ref{prob-oderRam-matchingK3} about numbers $R_<(K^<_3, \dots, K^<_3, M^<;q)$, where $M^<$ is an ordered matching on $n$ vertices and $q \geq 3$ is an integer.
Using a result by Alon and R\"{o}dl~\cite{alon2005sharp}, which gave a solution to a longstanding conjecture of Erd\H{o}s and Sós~\cite{erdos1980}, the known bounds are between $\Omega((n/ \log n)^{4/3})$ and $O (n^{q} \text{ poly} \log n)$.
These bounds are far apart and it is not clear where the truth should lie.

\begin{problem}[\cite{balGrin24,grin24}] Improve the upper and lower bounds on off-diagonal multicolor ordered Ramsey numbers $R_<(K^<_3, \dots, K^<_3, M^<;q)$.
\end{problem}

In the similar spirit of the conjecture by Erd\H{o}s and Sós~\cite{erdos1980}, Balko and Grinerov\'{a}~\cite{balGrin24,grin24} asked the following question. 

\begin{problem}[\cite{balGrin24,grin24}]
Are there ordered matchings $M^<$ on $n$ vertices such that $$\lim_{n \to \infty} \frac{R_<(K^<_3, K^<_3, M^<)}{R_<(K^<_3, M^<)} = \infty?$$
\end{problem}

For ordered graphs with bounded bandwidth, Gir\~{a}o, Janzer, and Janzer~\cite{gjj24} proved that $R_<(P^<_{k,n};q) \leq Dn^{Cq\log{q}}$ for some constants $D=D(k,q)>0$ and $C>0$.
They also mention some interesting consequences for the problem of finding monochromatic directed paths in multicolored tournaments and conjecture that their upper bound can be improved.

\begin{conjecture}[\cite{gjj24}]
For all positive integers $k,n,q$, there are constants $C>0$ and $D=D(k,q)>0$ such that $R_<(P^<_{k,n};q) \leq Dn^{Cq}$.
\end{conjecture}

The multicolor setting was also considered by Neidinger and West~\cite{neiWes19}, who provided various lower bounds on ordered Ramsey numbers of so-called \emph{2-ichromatic ordered graphs}, and by Balko, Cibulka, Kr\'{a}l, and Kyn\v{c}l~\cite{bckk13}, who showed the following dichotomy: $R_<(G^<; q)$ is either polynomial in $q$ or exponential in $q$ for every ordered graph $G^<$.
Bar\'at, Gy\'arf\'as, and T\'oth~\cite{barGyaTot24} also considered multicolor ordered Ramsey numbers of various classes of ordered mathchings and their connections to graph drawings.
However, overall, the $q$-color ordered Ramsey numbers remain quite unexplored.

\section{Hypergraph Ordered Ramsey Numbers}
\label{sec-orderedHypergraphs}

In full generality, Ramsey's theorem applies not only for graphs but also for \emph{$k$-uniform hypergraphs}, that is, pairs $(V, E)$ where $V$ is a finite set of \emph{vertices} and $E$ is a set of $k$-tuples from $V$ called \emph{edges}.
We use $K^{(k)}_n$ to dentote the complete $k$-uniform hypergraph $(V,\binom{V}{k})$.
Ramsey's theorem then states that, for all integers $k,q \geq 2$ and every $k$-uniform hypergraph $H$, there is a positive integer $N$ such that in every $q$-coloring of $K^{(k)}_N$ there a monochromatic copy of~$H$ forming a subhypergraph of $K^{(k)}_N$.
We denote the smallest such integer $N$ by $R(H;q)$ and call it the \emph{Ramsey number of $H$}.
If $q=2$, we write $R(H)$ instead of $R(H;2)$.

For $k \geq 3$, Ramsey numbers $R(K^{(k)}_n)$ are even less understood than their graph counterparts.
For example, it is only known that 
\begin{equation}
\label{eq-intro-hypergraphRamsey}
2^{\Omega(n^2)} \leq R(K^{(3)}_n) \leq 2^{2^{O(n)}},
\end{equation}
as shown by Erd\H{o}s, Hajnal, and Rado~\cite{erdosRado65}.
A famous conjecture of Erd\H{o}s, for whose proof Erd\H{o}s offered \$500 reward, states that there is a constant $c>0$ such that $R(K^{(3)}_n) \geq 2^{2^{cn}}$.
The case $k=3$ is of particular importance as if one determines the growth rate of $R(K^{(3)}_n)$ precisely, then the so-called \emph{Stepping-up lemma} by Erd\H{o}s and Rado~\cite{erdosRado52} would determine the growth rate of $R(K^{(k)}_n)$ for every $k \geq 4$.

The \emph{tower function} $t_h(x)$ of \emph{height} $h-1$ is defined by the recursive formula $t_1(x)=x$ and $t_h(x) =  2^{t_{h-1}(x)}$ for every $h \geq 2$.
Using~\eqref{eq-intro-hypergraphRamsey} together with the Stepping-up lemma, one can derive lower bounds on $R(K^{(k)}_n)$ for $k \geq 4$.
However, there is a difference of one exponential between known upper and lower bounds.
More precisely, we have
\begin{equation}
\label{eq-intro-tower}
t_{k-1}(\Omega(n^2)) \leq R(K^{(k)}_n) \leq t_k(O(n))
\end{equation}
for every fixed $k \geq 3$.

Perhaps surprisingly, the Ramsey number $R(H)$ of every $k$-uniform hypergraph $H$ with bounded $k$ and with bounded maximum degree is at most linear in the number of vertices of~$H$~\cite{crst83,cfs09,cnko08,cnko09,ishingami07,nsrs08}.

The \emph{ordered $k$-uniform hypergraph $H^<$} is a pair $(H,<)$ consisting of a $k$-uniform hypergraph $H$ and a linear ordering $<$ of its vertex set.
The notion of an \emph{ordered subhypergraph} is analogous to its graph counterpart.
Again, there is a unique ordered complete $k$-uniform hypergraph on $n$ vertices up to isomorphism and we denote it by $K^{<(k)}_n$.

The \emph{ordered Ramsey number} $R_<(H^<,G^<)$ of two ordered $k$-uniform hypergraphs $H^<$ and $G^<$ is the smallest $N \in \mathbb{N}$ such that every red-blue coloring of the hyperedges of $K^{<(k)}_N$ contains a blue copy of $H^<$ or a red copy of~$G^<$ forming an ordered subhypergraph.
In the \emph{diagonal case} $H^< = G^<$, we just write $R_<(H^<)$ instead of $R_<(H^<,H^<)$.
We can also naturally extend the definition to $q \geq 2$ colors, in which case we use $R_<(H^<;q)$ to denote the $q$-color variant of $R_<(H^<)$.

Similarly as for ordered graphs, the ordered Ramsey number of every $n$-vertex $k$-uniform hypergraph is bounded from above by $R(K^{(k)}_n)$.
In particular, these numbers are always finite by~\eqref{eq-intro-tower}.
It is also easy to see that the ordered Ramsey numbers of $k$-uniform hypergraphs grow at least as fast as the standard Ramsey numbers.

\subsection{Monotone Paths and the Erd\H{o}s--Szekeres Theorem}
\label{subsec-monotonePaths}

We start with one example, the original motivation for a systematic study of ordered Ramsey numbers.
We showed in Section~\ref{sec-orderedGraphs} that the Erd\H{o}s--Szekeres lemma is a consequence of a stronger Ramsey statement about monotone paths.
There is a similar connection between the celebrated \emph{Erd\H{o}s--Szekeres theorem} and ordered Ramsey numbers of 3-uniform hypergraphs.

To state the Erd\H{o}s--Szekeres theorem, a foundational result in discrete geometry, we first need to introduce some definitions.
A finite set $P$ of points in the plane is in \emph{general position} if no three points from $P$ lie on a common line.
A finite set of points is in \emph{convex position} if its points form vertices of a convex polygon.

\begin{theorem}[The Erd\H{o}s--Szekeres theorem~\cite{erdosSzekeres35}]
\label{thm-ES}
For every positive integer $n$, there is a positive integer $N$ such that every set of at least $N$ points in general position in the plane contains $n$ points in convex position.
\end{theorem}

We use $ES(n)$ to denote the smallest such integer $N$.
The statement of the Erd\H{o}s--Szekeres theorem is a generalization of Esther Klein’s problem, which was named the \emph{Happy Ending Problem} by Paul Erd\H{o}s, as it eventually led to the marriage of George Szekeres and Esther Klein.
Erd\H{o}s and Szekeres provided two proofs of this famous result.
One is an application of Ramsey's theorem and yields a rather poor upper bound on the function $ES(n)$.
The other proof uses more geometry and gives the estimate 
\begin{equation}
\label{eq-ES}
ES(n) \leq \binom{2n-4}{n-2}+1
\end{equation}
for every $n \geq 2$.
Already in 1935, Erd\H{o}s and Szekeres believed that this bound could be significantly improved.
Based on their results for $n=2,3,4$, they posed the famous and still open \emph{Erd\H{o}s--Szekeres conjecture}, for whose proof Erd\H{o}s offered \$500 reward.

\begin{conjecture}[The Erd\H{o}s--Szekeres conjecture~\cite{erdosSzekeres35}]
\label{conj-intro-ES}
For every integer $n \geq 2$, we have
\[ES(n) = 2^{n-2}+1.\]
\end{conjecture}

In the 1960s, Erd\H{o}s and Szekeres~\cite{erdosSzekeres60} supported Conjecture~\ref{conj-intro-ES} with the lower bound 
\begin{equation}
\label{eq-intro-ES-lower}
ES(n) \geq 2^{n-2}+1.
\end{equation}
Despite several attempts over the years~\cite{chung98,klei98,mojVla16,norYud16,petersSzekeres06,totVal98}, it is still open to decide whether the upper bound $ES(n) \leq 2^{n-2}+1$ holds.
It is only known that the conjecture is true for $n \leq 6$~\cite{petersSzekeres06}.
However, there was a breakthrough by Suk~\cite{suk17} in 2017, who proved a very close estimate $ES(n) \leq 2^{n+o(n)}$.

As mentioned, the Erd\H{o}s--Szekeres theorem, more precisely the bound~\eqref{eq-ES}, follows from known estimates on ordered Ramsey numbers of monotone 3-uniform paths. 
For an integer $k \geq 2$, the \emph{monotone $k$-uniform path} on $n$ vertices, denoted by $MP^{<(k)}_n$, is an ordered $k$-uniform $n$-vertex hypergraph with edges formed by $k$-tuples of consecutive vertices in~$<$; see part~(b) of Figure~\ref{fig-monPath}.
Note that the monotone path $MP^<_n$ corresponds to $MP^{<(2)}_n$.
The monotone paths are also called \emph{tight paths} in the literature.

We can again rather easily show that $ES(n) \leq R_<(MP^{<(3)}_n)$ by coloring triples of points according to their orientations with the ordering given by $x$-coordinates; see~\cite{bckk13,clfs17}.
The bound~\eqref{eq-ES} then follows from the fact 
\begin{equation}
\label{eq-monotone3UnifPaths}
R_<(MP^{<(3)}_n) = \binom{2n-4}{n-2}+1
\end{equation}
for every $n \geq 2$, which was proved by Moshkovitz and Shapira~\cite{moshkovitz12}.
Moreover, several other interesting geometric applications of estimates on $R_<(MP^{<(k)}_n)$ for $k \geq 3$ appeared, for example, variants of the Erd\H{o}s--Szekeres Theorem for convex bodies~\cite{fox12}, higher-order Erd\H{o}s--Szekeres theorems~\cite{elias11}, and Ramsey-type results~\cite{balko19} for \emph{$k$-intersecting pseudoconfigurations of points}~\cite{miya17}, \emph{$C_d$-arrangements of $n$ pseudohyperplanes} in $\mathbb{R}^{d}$~\cite{felWei01}, and extensions of the cyclic arrangement of hyperplanes with a pseudohyperplane~\cite{zieg93}.
Formulations of the Erd\H{o}s--Szekeres conjecture in terms of hypergraph ordered Ramsey numbers were also used in computer attacks on this conjecture and its generalizations~\cite{baeBal24,balVal17,petersSzekeres06}.

Given this motivation, the ordered Ramsey numbers $R_<(MP^{<(k)}_n)$ have been quite intensively studied~\cite{balko19,choudum02,elias11,fox12,milans12,moshkovitz12} even for higher uniformities and their growth rate is nowadays well understood.
Moshkovitz and Shapira~\cite{moshkovitz12} showed the following strong result.

\begin{theorem}[\cite{moshkovitz12}]
\label{thm-moshkovitzShapira}
For all positive integers $n$ and $k \geq 3$, 
\[
R_<(MP^{<(k)}_{n+k-1}) = t_{k-1}((2-o(1))n),
\]
where $t_{k-1}(\cdot)$ is the tower function of height $k-2$.
\end{theorem}

Note that the tower is of height one smaller than the tower-type upper bound on $R(K^{(k)}_n)$ from~\eqref{eq-intro-tower} obtained with the Stepping-up lemma.
In fact, Moshkovitz and Shapira~\cite{moshkovitz12} proved
\[R_<(MP^{<(k)}_{n+k-1})=\rho_k(n)+1,\]
where $\rho_k(n)$ is the number of \emph{line partitions of $n$ of order $k$} (see~\cite{moshkovitz12} for definitions).
For $k=3$, this gives the exact formula $R_<(MP^{<(3)}_n) = \binom{2n-4}{n-2}+1$.

Moshkovitz and Shapira~\cite{moshkovitz12} note that $R_<(MP^{<(3)}_n;q)-1$ is the number of antichains in $[n-2]^q$. 
In particular, for $n=4$, we get that $R_<(MP^{<(3)}_4;q)-1$ equals the \emph{Dedekind numbers}, the numbers of antichains in the Boolean poset $[2]^q$.
The problem of estimating the number of antichains in $[n]^q$ is a natural extension
of the classical problem of estimating the Dedekind numbers.
This problem from enumerative combinatorics was recently explored by Falgas-Ravry, R\"{a}ty, and Tomon~\cite{frt23} who resolved a conjecture of Moshkovitz and Shapira~\cite{moshkovitz12} about the number of antichains in $[n]^q$ and, consequently, about the growth rate of $R_<(MP^{<(3)}_n;q)$.
Their bound is very strong and even determines the leading constant in the exponent.

\begin{theorem}[\cite{frt23}]
For every $\varepsilon > 0$ if $q$ is a sufficiently large integer, then for every integer $n \geq 2$
\[(1-\varepsilon) \sqrt{\frac{6}{\pi(n^2-1)q}} \cdot n^q\leq \log{R_<(MP^{<(3)}_{n+2};q)} \leq (1+\varepsilon) \sqrt{\frac{6}{\pi(n^2-1)q}} \cdot n^q.\]
\end{theorem}

We conclude this subsection with a problem about $R_<(MP^{<(4)}_n)$ mentioned by Moshkovitz and Shapira~\cite{moshkovitz12}.
While we know the exact formula for $R_<(MP^{<(3)}_n)$ by~\eqref{eq-monotone3UnifPaths}, Moshkovitz and Shapira~\cite{moshkovitz12} argue that it might be hard to determine the asymptotic of the logarithm of $R_<(MP^{<(4)}_n)$.
This problem is connected to the extremal theory of partially ordered sets as $R_<(MP^{<(4)}_n)-1$ is exactly the number of antichains in the restricted \emph{Young’s lattice}~\cite{moshkovitz12}.

\subsection{Connections to Hypergraph Ramsey Numbers}

Here, we mention two interesting connections between ordered Ramsey numbers and hypergraph Ramsey numbers. 
A surprising relation between off-diagonal ordered Ramsey numbers of monotone $k$-uniform paths and the classical Ramsey numbers $R(K^{<(k-1)}_n)$ was discovered by Mubayi and Suk~\cite{mubSuk17}.
They showed the following estimate, which in the case $q=2$, $k = O(1)$, and $n$ tending to infinity shows that estimating $R_<(K^{<(k)}_n,MP^{<(k)}_{k+1})$ is equivalent to the famous problem of determining the tower growth rate of the Ramsey numbers $R(K^{(k-1)}_n)$ of complete $(k-1)$-uniform hypergraphs~\eqref{eq-intro-tower}.

\begin{theorem}[\cite{mubSuk17}]
For all integers $n > k \geq 2$,
\[R(K^{(k-1)}_{\lfloor n/q \rfloor};q) \leq R_<(K^{<(k)}_n,MP^{<(k)}_{k-q+1}) \leq  R(K^{(k-1)}_n;q).\]
\end{theorem}

Mubayi and Suk~\cite{mubSuk17} also introduced the following strengthening of a conjecture of Erd\H{o}s and Hajnal~\cite{erdRado72} about off-diagonal Ramsey numbers of complete $k$-uniform hypergraphs.

\begin{conjecture}[\cite{mubSuk17}]
For any fixed integer $k \geq 4$, \[R_<(K^{<(k)},MP^{<(k)}_{k+1}) \geq t_{k-1}(\Omega(n)).\]
\end{conjecture}

For more problems and results about off-diagonal ordered Ramsey numbers of monotone $k$-uniform paths, we recommend a recent survey by Mubayi and Suk~\cite{mubSuk20}; see also~\cite{gishJinSud24,mubSuk24} for some newer results.

Another connection was discovered by Conlon, Fox, Lee, and Sudakov~\cite{clfs17} who showed that for any 3-uniform hypergraph $H$, there is a family of ordered graphs $\mathcal{S}_H$ such that the Ramsey number of $H$ is
bounded in terms of the ordered Ramsey number of the family $\mathcal{S}_H$.
Here, the \emph{ordered Ramsey number $R_<(\mathcal{F})$ of a family $\mathcal{F}$} of ordered graphs is the smallest positive integer $N$ such that every 2-coloring of the edges of $K^<_N$ contains a monochromatic ordered copy of some ordered graph from $\mathcal{F}$.

For an ordered graph $G^<$ with the vertex set $[n]$, let $T(G^<)$ be a 3-uniform hypergraph on vertex set $[n+1]$ obtained by taking all triples whose first pair is an edge of $G^<$.
For a 3-uniform hypergraph $H$ on $n + 1$ vertices, we let $\mathcal{S}_H$ be the collection of ordered graphs $G^<$ on $[n]$ such that $H$ is a subhypergraph of $T(G^<)$.
Conlon, Fox, Lee, and Sudakov~\cite{clfs17} then related upper bounds on Ramsey numbers of 3-uniform hypergraphs to ordered Ramsey numbers by proving the following result.

\begin{theorem}[\cite{clfs17}]
\label{thm-ordRamHyper-connection}
Every 3-uniform hypergraph $H$ satisfies
\[R(H) \leq 2^{\binom{R_<(\mathcal{S}_H)}{2}} + 1.\]
\end{theorem}

Conlon, Fox, Lee, and Sudakov~\cite{clfs17} expect that the bound from Theorem~\ref{thm-ordRamHyper-connection} is close to sharp in many cases.
For example, the choice $H = K^{(3)}_{n+1}$ satisfies $\mathcal{S}_H = \{K^<_n\}$, for which Theorem~\ref{thm-ordRamHyper-connection} produces the double-exponential bound from~\eqref{eq-intro-hypergraphRamsey}, which is believed to be tight.
However, Conlon, Fox, Lee, and Sudakov~\cite{clfs17} constructed some cases where the bound from Theorem~\ref{thm-ordRamHyper-connection} is far from the truth.

\subsection{Bounded degrees and interval chromatic number}

Very little is known about ordered Ramsey numbers of general ordered $k$-uniform hypergraphs with $k \geq 3$ as ordered Ramsey numbers have been studied mostly for ordered graphs only.
We primarily focus on the case $k=3$ here.

A simple probabilistic argument provides the lower bound $R_<(H^<) \geq 2^{\Omega(n^2)}$ for every ordered 3-uniform hypergraph with $n$ vertices and $\Omega(n^3)$ hyperedges, which is of the same asymptotic growth rate as we have in the lower bound on $R_<(K^{<(3)}_n)$ by~\eqref{eq-intro-hypergraphRamsey}.
Thus, we consider mostly sparse 3-uniform hypergraphs.

The \emph{degree} of a vertex $v$ in a hypergraph $H$ is the number of hyperedges of $H$ that contain~$v$.
It follows from a result by Moshkovitz and Shapira~\cite{moshkovitz12} that there are ordered 3-uniform hypergraphs $H^<$ on $n$ vertices with maximum degree 3 such that $R_<(H^<) \geq 2^{\Omega(n)}$.

Therefore, to obtain smaller upper bounds on the ordered Ramsey numbers, it is necessary to bound other parameters besides the maximum degree.
A natural choice is again the interval chromatic number, which is defined analogously as for ordered graphs.
Recall that for ordered graphs, bounding both parameters indeed helps, as the ordered Ramsey number $R_<(G^<)$ of every $n$-vertex ordered graph $G^<$ with bounded maximum degree $d$ and bounded interval chromatic number $\chi$ is at most polynomial in the number of vertices by~\eqref{eq-ordRam-degInt}.
Theorem~\ref{thm-oderRam-degreeIntNum2} actually gives a stronger estimate
\[
R_<(G^<,K^<_\chi(n)) \leq n^{32 d \log{\chi}}.
\]

A natural question is whether we can get similar bounds for ordered $k$-uniform hypergraphs with $k \geq 3$.
For integers $k \geq 2$ and $\chi \geq k$, we use $K^{(k)}_\chi(n)$ to denote the \emph{complete $k$-uniform $\chi$-partite hypergraph}, that is, the vertex set of $K^{(k)}_\chi(n)$ is partitioned into $\chi$ sets of size $n$ and every $k$-tuple with at most one vertex in each of these parts forms a hyperedge.
Let $K^{<(k)}_\chi(n)$ be the ordering of~$K^{(k)}_\chi(n)$ in which the color classes form consecutive intervals.
Conlon, Fox, and Sudakov~\cite{cfs11} showed that, for all positive integers $\chi \geq 3$ and $n$, \[R(K^{(3)}_\chi(n)) \leq 2^{2^{2R}n^2},\]
where $R = R(K_{\chi-1})$.
Since every ordering of $K^{(3)}_\chi(\chi n)$ contains an ordered subhypergraph isomorphic to $K^{<(3)}_\chi(n)$ and every ordered $3$-uniform hypergraph on $n$ vertices with interval chromatic number $\chi$ is an ordered subhypergraph of $K^{<(3)}_\chi(n)$, we obtain 
$R_<(H^<) \leq 2^{2^{2R}\chi^2n^2}$,
for every ordered $3$-uniform hypergraph $H^<$ on $n$ vertices with interval chromatic number $\chi$ where $R = R(K_{\chi-1})$.
In particular, if the interval chromatic number $\chi$ of $H^<$ is fixed, we have $
R_<(H^<) \leq 2^{O(n^2)}$.

Note that the last bound is asymptotically tight for dense ordered hypergraphs with bounded interval chromatic number.
Thus, Balko and Vizer~\cite{balViz21} considered the sparse case in which we additionally bound the maximum degree, similarly as in the setting from Subsection~\ref{subsec-boundedDegInt}. 
More specifically, the first nontrivial case for ordered $3$-uniform hypergraphs with interval chromatic number $3$, for which they obtained an estimate with a subquadratic exponent.

\begin{theorem}[\cite{balViz21}]
\label{thm-ordRamHyper-3UnifMaxDegIntChr}
Let $H^<$ be an ordered 3-uniform hypergraph on $n$ vertices with maximum degree $d$.
Then there are constants $C=C(d)$ and $c>0$ such that
\begin{equation}
\label{eq-ordRamHyper-diagonal}
R_<(H^<,K^{<(3)}_3(n)) \leq 2^{O(n^{2-1/(1+cd^2)})}.
\end{equation}
\end{theorem}

The main idea of the proof of Theorem~\ref{thm-ordRamHyper-3UnifMaxDegIntChr} is based on an embedding lemma from~\cite{cfs12}, where the authors study Erd\H{o}s--Hajnal-type theorems for $3$-uniform tripartite hypergraphs.
Theorem~\ref{thm-ordRamHyper-3UnifMaxDegIntChr} immediately gives $R_<(H^<) \leq 2^{O(n^{2-\varepsilon})}$ for any ordered 3-uniform hypergraph $H^<$ on $n$ vertices with maximum degree $d$ and with interval chromatic number $3$.

The upper bound~\eqref{eq-ordRamHyper-diagonal} is quite close to the truth, as even when $H^<$ is fixed we get a superexponential lower bound on $R_<(H^<,K^{<(3)}_3(n))$, as shown by Fox and He~\cite{foxhe19} and independently by Balko and Vizer~\cite{balViz21}. 

\begin{theorem}[\cite{balViz21,foxhe19}]
\label{thm-ordRamHyper-3UnifMaxDegLower}
For every $t \geq 3$ and every positive integer $n$, we have
\[R_<(K^{<(3)}_{t+1},K^{<(3)}_3(n)) \geq 2^{\Omega(n \log{n})}.\]
\end{theorem}

For ordered hypergraphs of uniformity $k>3$, we recall that it follows from~\eqref{eq-intro-tower} that their ordered Ramsey numbers can grow as a tower of height $k-1$.
By modifying a result of Conlon, Fox, and Sudakov~\cite{cfs10}, Balko and Vizer~\cite{balViz21} showed that we do not have a tower-type growth rate for $R_<(H^<)$ once the uniformity and the interval chromatic number of $H^<$ are bounded.

\begin{proposition}[\cite{balViz21}]
\label{prop-ordRamHyper-generalBound}
Let $\chi,k$ be fixed integers with $\chi \geq k \geq 2$ and let $H^<$ be an ordered $k$-uniform hypergraph on $n$ vertices with interval chromatic number~$\chi$.
Then,
\[
R_<(H^<) \leq 2^{O(n^{\chi-1})}.
\]
\end{proposition}

The ordered Ramsey numbers for ordered $k$-uniform hypergraphs with $k \geq 3$ and the edge-ordered Ramsey numbers are quite unexplored, so there is a plenty of open problems.

Theorems~\ref{thm-ordRamHyper-3UnifMaxDegIntChr} and~\ref{thm-ordRamHyper-3UnifMaxDegLower} give estimates on the ordered $R_<(G^<, K^{<(3)}_3(n))$ and although the exponents in the bounds are reasonably close, there is still a gap between them and it would be interesting to close it.

\begin{problem}[\cite{balViz21}]
Let $d$ be a fixed positive integer and let $H^<$ be an ordered $3$-uniform hypergraph on $n$ vertices with maximum degree $d$.
Close the gap between the lower and upper bounds on $R_<(H^<, K^{<(3)}_3(n))$.
\end{problem}

Another interesting problem is to extend the upper bound with subquadratic exponent from Theorem~\ref{thm-ordRamHyper-3UnifMaxDegIntChr} to ordered $3$-uniform hypergraphs with bounded maximum degree and fixed interval chromatic number that is larger than $3$.

\begin{problem}[\cite{balViz21}]
Let $d$ and $\chi$ be fixed positive integers.
Is there an $\varepsilon = \varepsilon(d,\chi)>0$ such that, for every ordered $3$-uniform hypergraph $H^<$ on $n$ vertices with maximum degree $d$ and with interval chromatic number $\chi$, we have
\[R_<(H^<) \leq 2^{O(n^{2-\varepsilon})}?\]
\end{problem}

In general, we are not aware of any nontrivial upper bounds on ordered Ramsey numbers of ordered $3$-uniform hypergraphs with bounded maximum degree. 

\begin{problem}[\cite{balViz21}]
What is the upper bound on ordered Ramsey numbers of ordered $3$-uniform hypergraphs with bounded maximum degree?
\end{problem}

\section{Edge-ordered Ramsey Numbers}
\label{sec-edgeOrdered}

Besides vertices, we can also order the edges of a given graph and then study its extremal properties.
An \emph{edge-ordered graph} $G^\prec$ is a pair $(G,\prec)$ consisting of a graph $G=(V, E)$ and a linear ordering $\prec$ of the set of edges~$E$.
We use the term \emph{edge-ordering of $G$} for the ordering $\prec$ and also for $G^\prec$.
An edge-ordered graph $G^\prec$ is an \emph{edge-ordered subgraph} of an edge-ordered graph $H^\prec$ if $G$ is a subgraph of $H$ and the edge-ordering of $G$ is a suborder of the edge-ordering of $H$.
We say that $G^\prec$ and $H^\prec$ are \emph{isomorphic} if there is a graph isomorphism between $G$ and $H$ that also preserves the edge-orderings.

Motivated by extremal results for edge-graphs obtained by Gerbner et al.~\cite{gmnptv19}, Balko and Vizer~\cite{balViz20} introduced Ramsey numbers for edge-ordered graphs.
The \emph{edge-ordered Ramsey number} $R_\prec(G^\prec)$ of an edge-ordered graph $G^\prec$ is the minimum positive integer $N$ such that there exists an edge-ordering $K^\prec_N$ of~$K_N$ such that every 2-coloring of the edges of~$K^\prec_N$ contains a monochromatic copy of $G^\prec$ as an edge-ordered subgraph of $K^\prec_N$.

Note that the definition of edge-ordered Ramsey numbers is defined quite differently than ordered Ramsey numbers as the edge-ordering of the complete graph whose edges are being colored depends on the given edge-ordered graphs.
This is necessary, as otherwise there might be an edge-ordered graph $G^\prec$ and an edge-order of $K^\prec_N$ such that $G^\prec$ is not an edge-ordered subgraph of $K^\prec_N$.

The finiteness of ordered Ramsey numbers was quite easy to show, as it followed from the finiteness of standard Ramsey numbers.
This is not the case for edge-ordered Ramsey numbers, where it takes some effort to prove that these numbers are finite for any pair of edge-ordered graphs.
This was proved by Balko and Vizer~\cite{balViz20}.

\begin{theorem}[\cite{balViz20}]
\label{thm-ordRamHyper-edgeOrderedRamseyFinite}
For every edge-ordered graph $G^\prec$, the edge-ordered Ramsey number $R_\prec(G^\prec)$ is finite.
\end{theorem}

Theorem~\ref{thm-ordRamHyper-edgeOrderedRamseyFinite} also follows from a recent deep result of Hubi\v{c}ka and Ne\v{s}et\v{r}il~\cite[Theorem~4.33]{hubNes19} about Ramsey numbers of general relational structures.
In comparison, the proof of Theorem~\ref{thm-ordRamHyper-edgeOrderedRamseyFinite} is less general, but it is much simpler and produces better and more explicit bound on $R_\prec(G^\prec)$.
The proof of Theorem~\ref{thm-ordRamHyper-edgeOrderedRamseyFinite} yields a stronger induced-type statement where additionally the ordering of the vertex set is fixed and the colorings can use an arbitrary number of colors.

However, the bound on the edge-ordered Ramsey numbers obtained in the proof of Theorem~\ref{thm-ordRamHyper-edgeOrderedRamseyFinite} is enormous, it grows faster than, for example, a tower function of any fixed height.
Fox and Li~\cite{foxLi20} improved the bound on edge-ordered Ramsey numbers to a single exponential type.

\begin{theorem}[\cite{foxLi20}]
\label{thm-ordRamHyper-foxLi}
For each positive integer $n$, there is an edge-ordered graph $G^\prec$ on $N = 2^{100n^2\log^2{n}}$ vertices such that, for every 2-coloring of the edges of $G^\prec$, there exists a monochromatic subgraph containing a copy of every $n$-vertex edge-ordered graph.
\end{theorem}

In particular, if $G^\prec$ is an edge-ordered graph on $n$ vertices, then
\[R_\prec(G^\prec) \leq 2^{100n^2\log^2{n}}.\]

Fox and Li~\cite{foxLi20} also extended their result to more colors and proved the following polynomial upper bound on edge-ordered Ramsey numbers of edge-ordered graphs of bounded degeneracy, improving earlier estimates by Balko and Vizer~\cite{balViz20}.

\begin{theorem}[\cite{foxLi20}]
\label{thm-ordRamHyper-foxLiDegen}
If $G^\prec$ is an edge-ordered $d$-degenerate graph on $n$ vertices, then
\[R_\prec(G^\prec) \leq n^{600d^3\log{(d+1)}}.\]
\end{theorem}

For edge-ordered graphs, a major open question is to find a better estimate for edge-ordered Ramsey numbers of edge-ordered complete graphs on $n$ vertices.
In particular, Fox and Li~\cite{foxLi20} asked whether the exponential lower bound is tight for general edge-orderings of $K_n$.
For sparser edge-ordered graphs, Fox and Li~\cite{foxLi20} conjectured that the upper bound from Theorem~\ref{thm-ordRamHyper-foxLiDegen} can be improved.

\begin{conjecture}[\cite{foxLi20}]
\label{prob-ordRamHyper-edgeDegen}
If $H^\prec$ is an edge-ordered $d$-degenerate graph on $n$ vertices, then $R_\prec(H^\prec) \leq n^{O(d)}$.
\end{conjecture}

Currently, there are no known edge-ordered $d$-degenerate graphs with superlinear edge-ordered Ramsey numbers.
However, Fox and Li~\cite{foxLi20} believe that there are such examples and conjectured that the upper bound from Conjecture~\ref{prob-ordRamHyper-edgeDegen} is tight up to the constant in the exponent.

The edge-ordered Ramsey numbers can be naturally extended to edge-ordered $k$-uniform hypergraphs with any $k \geq 2$.
The existence of such numbers follows from a result of Hubi\v{c}ka and Ne\v{s}et\v{r}il~\cite{hubNes19}, but the resulting bounds are enormous.
Fox and Li~\cite{foxLi20} posed a natural problem to give a better estimate on such hypergraph edge-ordered Ramsey numbers.

\bibliography{bibliography}
\bibliographystyle{plainnat}

\end{document}